\documentclass[10pt, reqno]{amsart}
\usepackage{amssymb,latexsym}
\usepackage{eucal}
\newtheorem{thm}{Theorem}
\newtheorem{lemma}{Lemma}
\newtheorem*{proposition*}{}

 \DeclareMathOperator{\Ker}{\textup{Ker}}
 
\DeclareMathOperator{\Aut}{\textup{Aut}}

\newcommand{\ph}{{\varphi}}

\newcommand{\X}{\mathcal{X}}
\newcommand{\Ss}{\mathcal{S}}
\newcommand{\PP}{\mathcal{P}}
\newcommand{\R}{\mathcal{R}}
\newcommand{\U}{\mathcal{U}}
\newcommand{\HH}{\mathcal{H}}
\newcommand{\B}{\mathcal{B}}
\newcommand{\Q}{\mathcal{Q}}
\newcommand{\T}{\mathcal{T}}
\newcommand{\A}{\mathcal{A}}
\newcommand{\C}{\mathcal{C}}
\newcommand{\Y}{\mathcal{Y}}
\newcommand{\K}{\mathcal{K}}
\newcommand{\G}{\mathcal{G}}
\newcommand{\D}{\mathcal{D}}

\newcommand{\V}{\mathcal{V}}
\newcommand{\F}{\mathcal{F}}

\newcommand{\p}{\partial}
\newcommand{\E}{\textup{E}}
\newcommand{\kp}{\tau}
\newcommand{\ra}{\rangle}
\newcommand{\la}{\langle}
\newcommand{\wtl}{\widetilde}

\begin{document}

\title[Embedding free Burnside groups]
{Embedding free Burnside groups in finitely presented groups}
\author{S.V. Ivanov}
\address{Department of Mathematics\\
University of Illinois \\
Urbana,   IL 61801, U.S.A.}
\email{ivanov@math.uiuc.edu}
\thanks{Supported in part by NSF grant   DMS 00-99612}
\begin{abstract}
We construct an embedding of a free Burnside group $B(m,n)$ of odd $n >
2^{48}$ and rank $m >1$ in a finitely presented group with some special
properties. The main application of this embedding is an easy construction of
finitely presented non-amenable groups without noncyclic free subgroups
(which provides a finitely presented counterexample to the von Neumann
problem on amenable groups). As another application, we construct weakly
finitely presented groups of odd exponent $n \gg 1$ which are not locally
finite.
\end{abstract}

\maketitle
\section{Introduction}

The class of amenable groups was introduced by von Neumann \cite{N29} to
study so-called Hausdorff-Banach-Tarski paradox  \cite{H14},  \cite{BT24}
(see also \cite{OS01}). There are several equivalent  definitions of group
amenability. For example, a group $G$ is called {\em amenable} (or
measurable) if  there is a measure $\mu$ on $G$ (that is, $\mu$ is a
nonnegative real function on subsets of $G$) such that $\mu(X \cup Y) =
\mu(X) + \mu(Y)$ whenever $X$, $Y$ are disjoint subsets of $G$, $\mu(G) =1$
and $\mu(gX) = \mu(X)$ for all $g \in G$, $X \subseteq G$. Clearly, finite
groups are amenable. von Neumann \cite{N29}  showed that the class of
amenable groups contains abelian groups, contains no free group $F_2$ of rank
2, and is closed under  operations of taking subgroups, extensions, and
(infinite) ascending unions. Day \cite{D57} and Specht  \cite{S60}  proved
that this class is closed under taking homomorphic images.

The problem, attributed to von Neumann, on the existence of non-amenable
groups without noncyclic free subgroups was solved in the affirmative by
Ol'shanskii \cite{O80}. Recall that the alternative "either amenable or
contains $F_2$", suggested by  the above von Neumann results, was proved by
Tits  \cite{T72} for linear groups over a field of characteristic 0. Adian
\cite{A82} showed that a  free $m$-generator Burnside group $B(m,n)$ of odd
exponent $n \ge 665$ is not amenable if $m >1$ (recall that $B(m,n)$ is the
quotient $F_m /F_m^n$, where $F_m$ is a free group of rank $m$; in
particular, $B(m,n)$ contains no free subgroups of rank $\ge 1$).

However, both Ol'shanskii's  and Adian's examples of non-amenable groups are
not finitely presentable and, in this connection, a "finitely presented"
version of the von Neumann problem (on the existence of finitely presented
non-amenable groups without noncyclic free subgroups) was promptly stated by
Cohen \cite{C82} and Grigorchuk \cite[Problem 8.7]{KN82}. To solve this
problem, Ol'shanskii and Sapir \cite{OS02} (see also \cite{OS01}) constructed
a finitely presented group $\G_{OS}$ such that $\G_{OS}$ is an extension of a
group  of odd exponent $n \gg 1$ by an infinite cyclic group (so $\G_{OS}$
has no noncyclic free subgroups) and $\G_{OS}$ contains a subgroup isomorphic
to a free 2-generator Burnside group $B(2,n)$ of exponent $n$ (and so, by
Adian's result \cite{A82},  $\G_{OS}$  is not amenable).  The construction of
the group $\G_{OS}$  is rather involved and lengthy; it is based upon a
modification of Ol'shanskii's method of graded diagrams \cite{O89} and uses
$S$-machines introduced by Sapir   \cite{SBR97}. In this regard, recall that
the first example of a finitely presented torsion-by-cyclic group $\G_G$ was
constructed by Grigorchuk \cite{G98}. However, the group $\G_G$ is amenable
(but not elementary amenable). The main goal of this article is to use the
techniques of HNN-extensions in the class of groups of odd exponent $n \gg
1$, as developed in \cite{I02b}, to give a significantly simpler, compared to
\cite{OS02}, construction of finitely presented non-amenable groups without
noncyclic free subgroups. Similar to examples of Grigorchuk and
Ol'shanskii-Sapir,  our finitely presented groups are also torsion-by-cyclic
(which guarantees the absence of noncyclic free subgroups).

Recall that a subgroup $H$ of a group $G$ is called a {\em $Q$-subgroup} if
for every normal subgroup $K$ of $H$ the normal closure $\la K \ra^G$ of $K$
in $G$ meets $H$ in $K$ (this property was introduced by Higman, B. Neumann,
and H. Neumann \cite{HNN49} to show that every countable group embeds in a
2-generator group). Here is our main technical result.

\begin{thm}
Let $n > 2^{48}$ be odd,  $m >1$, and let $ \G_U = \langle \U \; \| \;
(\ref{R1})-(\ref{R7}) \rangle $ be the  finitely  presented group whose
generating set is the alphabet $\U$, given by $(\ref{U})$, and whose defining
relations are $(\ref{R1})-(\ref{R7})$. Then the subgroup $\G_A = \la a_1,
\dots, a_m \ra$ of $\G_U$, generated by the images of $a_1, \dots$, $a_m$, is
a $Q$-subgroup of $\G_U$ naturally isomorphic to a free $m$-generator
Burnside group $B(m,n)$ of exponent $n$. Furthermore, the group $\G_A  = \la
a_1, \dots, a_m \ra$ naturally embeds in the quotient $\G_U/\G_U^n$ as a
$Q$-subgroup.
\end{thm}

Recall that, by the Higman embedding theorem \cite{H61}, a finitely generated
group $G$ embeds in a finitely presented group if and only if $G$ can be
recursively presented. This, in particular, means that a free $m$-generator
Burnside group $B(m,n)$ of exponent $n$, being recursively presented for all
$m$, $n$, embeds in a finitely generated group and the significance of
Theorem~1 is in the additional properties of this embedding.

As a corollary of Theorem 1 and results of article \cite{I02c} on
$Q$-subgroups of $B(m,n)$, we have the following theorem that strengthens the
main result of \cite{OS02}.
\begin{thm}
Let $n > 10^{78}$ be odd and the finitely presented group $\G_{T2}$ be
defined by presentation $(\ref{GT2})$. Then $\G_{T2}$ is an extension of a
group $\K_{T2}$ of exponent $n$  by an infinite cyclic group and $\K_{T2}$
contains a subgroup isomorphic to a free $2$-generator Burnside group
$B(2,n)$ of exponent $n$. In particular, $\G_{T2}$ is not amenable and
contains no noncyclic free subgroups. In addition, one can assume that if
$\V(B(2,n))$ is a given nontrivial verbal subgroup of $B(2,n)$ then
$\V(\K_{T2}) = \K_{T2}$.
\end{thm}

A group presentation $G = \la \C \; \| \; R=1, R \in \R \ra$ is called {\em
weakly finite} (and the group $G$ is termed {\em weakly finitely presented})
if every finitely generated subgroup $H$ of $G$, generated  by the images of
some words in $\C^{\pm 1}$, is naturally isomorphic to the subgroup of a
group $\G_{0} = \la \C_0 \, \| \,  R=1, R \in \R_0 \ra$, where $\C_0
\subseteq \C$, $\R_0 \subseteq \R$ are finite subsets depending upon $H$,
generated by  the images of the same words. For example, if the alphabet $\C$
is finite then $G = \la \C \, \| \, R=1, R \in \R \ra$  is  weakly finite if
and only if there is a finite subset $\R_1 \subseteq \mathcal R$ such that
$\la \R_1 \ra^{F(\C)} = \la \R \ra^{F(\C)}$, where $F(\C)$ is the free group
in the alphabet $\C$ and $\la \R_1 \ra^{F(\C)}$, $\la \R \ra^{F(\C)}$ are
normal closures in $F(\C)$  of $\R_1$, $\R$, respectively. It is also easy to
see that a group presentation $G = \la \C \, \| \, R=1, R \in \R \ra$ is
weakly finite if every finitely generated subgroup of $G$ is finitely
presentable (i.e., has a finite presentation; however, the converse does not
hold as can be seen from Theorem 3). In particular, if $G$ is a locally
finite group then $G$ is weakly finitely presented. It is rather natural to
ask whether there are periodic weakly finitely presented groups that are not
locally finite. This problem can be regarded  as a weakened version of a
long-standing problem, attributed to P.S.~Novikov, on the existence of
finitely presented infinite periodic groups.

The article \cite{I02a} presents a construction  of weakly finitely presented
periodic groups which are not locally finite thus making a progress towards
the Novikov problem. However, the construction of \cite{I02a} does not allow
to have  orders of elements of the group bounded and a construction of weakly
finitely presented groups of bounded exponent that are not locally finite
seems to be the next reasonable goal. We also remark that the construction of
\cite{I02a}  is based upon a technical embedding result which could be stated
exactly as Theorem 2  with the modification that the quotient $\G_U/\G_U^n$
is $\G_U/\G_U^{n^3}$ in  \cite{I02a} (in \cite{I02a} the group $\G_U$ is
given by simpler than (\ref{R1})--(\ref{R7}) relations and the estimate $n
> 2^{48}$ is $n > 2^{16}$ in \cite{I02a}). Thus,  Theorem 1 strengthens
the main construction of \cite{I02a} and it can be used (in a  similar to
\cite{I02a} fashion) to construct weakly finitely presented non-locally
finite groups of bounded exponent $n \gg 1$ (the question on the existence of
such groups was raised in \cite{I02a}).
\begin{thm}
Let $n > 10^{78}$ be odd. Then there is a weakly finitely presented group
$\G_{T3}$ of exponent $n$ which contains a subgroup isomorphic to $B(2,n)$.
\end{thm}

\section{Five Lemmas}

Let $\A = \{ a_1, \dots, a_m \}$, $m>1$, be an alphabet. We also consider
alphabets
\begin{gather*}
\B = \{ b_1, \dots, b_m \} ,  \ \  \bar \B = \{ \bar  b_1, \dots, \bar b_m \} , \  \
\D = \{ d_1, \dots, d_m \} ,  \ \ \bar \D = \{\bar  d_1, \dots, \bar  d_m \} , \\
\X = \{ x_1, \dots, x_m \} ,  \ \ \bar \X = \{\bar  x_1, \dots, \bar  x_m \} , \  \ \Y =
\{ y_1, \dots, y_m \} ,  \ \ \bar \Y = \{\bar  y_1, \dots, \bar  y_m \} ,
\end{gather*}
and  $\{ c, \bar c, e , \bar e, z \}$. Let $\U$ be the union of all these
pairwise disjoint alphabets
\begin{equation} \label{U}
\U = \A \cup \B \cup \bar \B \cup \D \cup \bar \D \cup \X \cup \bar \X \cup
\Y \cup \bar \Y \cup \{ c, \bar c, e, \bar e, z \} .
\end{equation}

Suppose that $\R_A$ is a set of words in $\A^{\pm 1} = \A \cup \A^{-1}$ and
\begin{equation}  \label{GA}
\G_A = \langle \A \; \| \; R=1, R \in \R_A \rangle
\end{equation}
is a group presentation whose set of generators is $\A$ and whose set of
defining words is $\R_A$. Consider the following defining relations in which
$[x, y] = xyx^{-1}y^{-1}$ is the commutator of $x, y$ and indices $i, j$ run
over $\{1, \dots, m\}$:
\begin{gather}
x_i c x_i^{-1} = cb_i, \ \ \bar x_i \bar c \bar x_i^{-1} = \bar c \bar b_i, \ \
y_i e y_i^{-1} = e d_i, \ \ \bar y_i \bar e \bar y_i^{-1} = \bar e \bar d_i,
\label{R1} \\
[x_i, b_j] = 1 , \ \  [\bar x_i, \bar b_j ] = 1 , \ \  [y_i, d_j] = 1, \ \
[\bar y_i, \bar d_j] = 1 , \label{R2} \\
c^n =1, \ \ \bar c^n = 1, \ \ e^n =1, \ \  \bar e^n = 1, \label{R3} \\
[b_i, \bar b_j ] = 1, \ \  [b_i, \bar c] = 1 , \ \ [\bar b_j, c] = 1,  \ \
[c, \bar c ] =
1 , \label{R4} \\
[d_i, \bar d_j] =1, \ \  [d_i, \bar e] =1 , \ \  [\bar d_j, e] = 1, \ \
[e, \bar e ] = 1
, \label{R5} \\
[a_i, d_j ]= 1, \ \  [a_i, \bar d_j]= 1, \ \ [a_i, y_j ]= 1, \ \
[a_i, \bar y_j]= 1, \ \
[a_i, e ]= 1, \ \  [a_i, \bar e]= 1, \label{R6} \\
z b_i \bar b_i z^{-1} = d_i \bar d_i a_i , \ \ z c \bar c z^{-1} = e \bar e
\label{R7} .
\end{gather}

Less formally,  these relations can be described as follows. First we
introduce relations $x_i c x_i^{-1} = cb_i$, $[x_i, b_j] = 1$, $i, j \in \{
1, \dots, m\}$, $c^n = 1$ (which copy relations of \cite{I02a}) on letters in
$\B \cup \X \cup \{c \}$. Then we consider "bar" copies of these letters and
relations in $\bar \B \cup \bar \X \cup \{\bar c \}$ and make letters in $\B
\cup \{c \}$ commute with letters in $ \bar \B \cup \{\bar c \}$. Then  we
introduce copies of all letters in $\B \cup  \bar \B \cup \X \cup  \bar \X
\cup \{c, \bar c \}$, denoted by $\D \cup  \bar \D \cup \Y \cup  \bar \Y \cup
\{e, \bar e \}$ (so that $\D$ is a copy of $\B$, $\bar Y$ is a copy of $\bar
X$ etc.), and make copies of already defined relations in $\B \cup  \bar \B
\cup \X \cup  \bar \X \cup \{c, \bar c \}$ to get their duplicates in $\D
\cup  \bar \D \cup \Y \cup  \bar \Y \cup \{e, \bar e \}$. Finally, we define
an alphabet $\A$, make letters in $\A$ commute with letters in $\D \cup  \bar
\D \cup \Y \cup  \bar \Y \cup \{e, \bar e \}$, and introduce one more letter
$z$ and $m+1$ more relations $z b_i \bar b_i z^{-1} = d_i \bar d_i a_i$,
$i=1, \dots, m$, $z c \bar c z^{-1} = d \bar d$.

Let a group $\G_U$ be given by a presentation whose alphabet is $\U$ and
whose set of defining relations consists of defining relations of
presentation (\ref{GA}) and relations (\ref{R1})--(\ref{R7}), thus
\begin{equation} \label{GU}
\G_U = \langle \U \; \| \; R = 1, \; R \in \R_A, \ (\ref{R1})-(\ref{R7}) \rangle .
\end{equation}

Denote the free group in the alphabet $\A$ by $F(\A)$ and let
$$
\alpha_A : F(\A) \to \G_A,  \quad  \alpha_U : F(\U) \to \G_U
$$
be natural homomorphisms.

A word $W$ in the alphabet $\A^{\pm 1}$ is called an $\A$-{\em word} and
denoted by $W = W(\A) = W(a_1, \dots, a_m)$. The substitution $a_i \to b_i$,
$i =1, \dots, m$, turns an $\A$-word into a $\B$-word $W(\B) = W(b_1, \dots,
b_m)$. Recall that if $G$ is a group then $G^n = \langle g^n \; | \; g \in G
\rangle$ denotes the subgroup of $G$ generated by all $n$th powers of
elements of $G$.

\begin{lemma} The map $a_i^{\alpha_A} \to a_i^{\alpha_U}$, $i =1, \dots, m$,
extends to a homomorphism $ \psi_0:  \G_A \to \G_U $ whose kernel $\Ker
\psi_0$ is $\G_A^n$.
\end{lemma}

{\em Proof.} This is analogous to the proof of Lemma A \cite{I02a}. First we
will show that $\G_A^n \subseteq \Ker \psi_0$. Let $W = W(\mathcal A)$ be an
$\A$-word. According to relations (\ref{R1})--(\ref{R3}), the  following
equalities hold in the group $\G_U$
$$
1 \overset  {\G_U}  = (W (\X) c W(\X)^{-1} )^n \overset {\G_U}   = ( c
W(\B))^n \overset {\G_U}  = (\bar c W(\bar \B))^n \overset {\G_U} = ( e
W(\D))^n \overset  {\G_U}  = (\bar e W(\bar \D))^n .
$$
Hence, by relations (\ref{R4}), we get
\begin{multline*}
1 \overset  {\G_U}  = ( c W(\B))^n (\bar c W(\bar \B))^n \overset {\G_U} = (
c \bar c W(\B)  W(\bar \B) )^n \overset {\G_U} =  ( c \bar c W(b_1 \bar b_1,
\dots, b_m \bar b_m) )^n .
\end{multline*}
Conjugating the last word by $z$, in view of relations (\ref{R7}),
(\ref{R6}), (\ref{R5}) we further have
\begin{multline*}
1 \overset  {\G_U}  = z ( c \bar c W(b_1 \bar b_1, \dots, b_m \bar b_m) )^n
z^{-1} \overset  {\G_U}  =  ( e \bar e W(d_1 \bar d_1 a_1,\dots, d_m \bar d_m
a_m) )^n \overset {\G_U}  = \\  \overset  {\G_U}  = ( e \bar e W(d_1 \bar
d_1, \dots, d_m \bar d_m) W(a_1, \dots, a_m) )^n  \overset  {\G_U}  = ( e W(
\D) \bar e W(\bar  \D) )^n W(\A)^n \overset  {\G_U}  =
\\  \overset  {\G_U}  = ( e W( \D))^n (\bar e W(\bar  \D) )^n W(\A)^n
\overset  {\G_U}  = W(\A)^n  .
\end{multline*}
Thus, $W(\A) \overset  {\G_U}  = 1$ for every $\A$-word $W(\A)$ and the
inclusion $\G_A^n \subseteq \Ker \psi_0$ is proven.

To prove the converse, consider a free Burnside group $B(\B \cup c, n)$ of
exponent $n$ which is freely generated by $\B \cup c$ ($B(\B \cup c, n)$ is
the quotient $F(\B \cup c)/F(\B \cup c)^n$). A  free Burnside group $B(\X,
n)$ is defined analogously. It is easy to see that one can consider a
semidirect product
$$
S(\B \cup c \cup \X) = B(\B \cup c, n) \leftthreetimes  B(\X, n)
$$
so that $x_i b_j x_i^{-1} = b_j$ , \ $x_i c x_i^{-1} = c b_i$,  \ $i, j \in
\{1, \dots, m\}$.

Semidirect products $S(\bar \B \cup \bar  c \cup  \bar \X)$ , \ $S(\D \cup e
\cup \Y)$, \ $S(\bar \D \cup \bar  e \cup  \bar \Y)$ are defined analogously.

Consider the direct products
\begin{equation} 
\label{PBcX} S(\B \cup c \cup \X) \times S(\bar \B \cup \bar  c \cup  \bar
\X) , \quad S(\D \cup e \cup \Y) \times S(\bar \D \cup \bar  e \cup  \bar \Y)
\times \G_A / \G_A^n
\end{equation}
and observe that the map $ b_i \bar b_i \to d_i \bar d_i a_i$,  \ $c \bar c
\to e \bar e$, \  $i =1, \dots, m$,  extends to an isomorphism of
corresponding subgroups of groups (\ref{PBcX}) (these subgroups are obviously
isomorphic to $B(\B \cup c, n)$). Therefore, one can consider an
HNN-extension $H_{z}$ of the free product of groups (\ref{PBcX}) with stable
letter $z$ whose relations are $ z b_i \bar b_i z^{-1} = d_i \bar d_i a_i$, \
$z c \bar c z^{-1} =  e \bar e$, \ $i =1, \dots, m$.

Now we observe that this HNN-extension $H_{z}$ has generating set $\U$ and
all relations (\ref{R1})--(\ref{R7}) hold in these generators. This means
that the natural homomorphism $F(\U) \to H_{z}$  naturally factors out
through $\G_U$. Since the group $\G_A/\G_A^n$ naturally embeds in $H_{z}$,
the group $\G_A/\G_A^n$ also embeds in $\G_U$ which implies the converse
inclusion $\Ker \psi_0 \subseteq \G_A^n$. Thus, $\Ker \psi_0 = \G_A^n$ and
Lemma 1 is proved. \qed
\smallskip

Let $P$ be a group of {\em exponent} $n$ (that is, $p^n =1$ for every $p \in
P$) and let $\T = \{ t_1, \dots, t_\ell \}$ be an alphabet, $\ell \geq 1$.
Let $P_{k,1}$, $P_{k,2}$ be two isomorphic subgroups of $P$, $k=1, \dots,
\ell$, and
$$
\rho_k : P_{k,1} \to P_{k,2}
$$
be a fixed isomorphism. Consider the following HNN-extension of $P$:
\begin{equation}
\label{HHT} \HH_T = \langle P, t_1, \dots, t_\ell \; \| \; t_k p t_k^{-1} =
\rho_k(p), p \in P_{k,1} , k = 1, \dots, \ell \rangle .
\end{equation}

Let $A \in \HH_T$, $\F(A)$ denote the maximal subgroup of $P \subset \HH_T$
that is normalized by $A$ and $\kp_A$ be the automorphism of $\F(A)$ which is
induced by conjugation by $A$. Consider the following property of $\HH_T$ in
which $N_E \ge 1$ is an integer parameter.

\begin{enumerate}
\item[(E)] For every $A \in \HH_T$ which is not conjugate in the group
$\HH_T$ given by (\ref{HHT}) to an element $p \in P$, the subgroup $\langle
\kp_A, \F(A) \rangle$ has exponent $n$ and equalities $A^{-k} p_0 A^{k} =
p_k$, where $p_k \in P$ and $k = 0, \dots, N_E$, imply that $p_k \in \F(A)$.
\end{enumerate}

It is proved in \cite{I02b} that if $n > 2^{48}$ is odd  and the property (E)
holds for $\HH_T$ with $N_E = [2^{-16}n]$, where $[2^{-16}n]$ is the integer
part of $2^{-16}n$, then the group $P$ naturally embeds in the quotient
$\HH_T / \HH_T^n$. This result and the techniques used in its proof will play
a central role in our subsequent arguments.
\smallskip

Let $B(\B, n)$, $B(\X, n)$ be free Burnside groups of exponent $n$ in alphabets
$\B$, $\X$, respectively, and $\PP_{B,X}$ be the direct product
$$
\PP_{B,X} = B(\B, n) \times B(\X, n) .
$$
Clearly, the map
$
x_i \to b_i x_i, \quad i =1, \dots, m,
$
extends to an isomorphism of corresponding subgroups of  $\PP_{B,X}$ and so we
can introduce an HNN-extension of $\PP_{B,X}$ with stable letter $c$:
\begin{equation}  \label{HHc}
\HH_c = \langle \PP_{B,X}, c  \; \| \; c^{-1} x_i c = x_i b_i, i = 1, \dots,
m \rangle .
\end{equation}

\begin{lemma} The group $\HH_c$, given by presentation $(\ref{HHc})$, has property $(E)$ with
$N_E =~3$.
\end{lemma}

{\em Proof.} Assume that $A \in \HH_c$, $A$ is not conjugate in $\HH_c$ to $p
\in \PP_{B,X}$ and $A^{-1}  p_{j-1} A = p_j$, where $p_{j-1}, p_j \in
\PP_{B,X}$ and $j=1,2,3$, $p_0 \neq 1$. Our goal is to show that $p_0 \in
\F(A)$ and the subgroup $\langle \kp_A, \F(A) \rangle$ has exponent $n$.

We will keep the terminology and notation of article \cite{I02b} (which is
based upon \cite{I94}, see also \cite{I92b}, \cite{I98}, \cite{I02a}).

Without loss of generality, we can assume that if $A = B$ in $\HH_c$ then
$|A|_c \leq |B|_c$,  where $|A|_c$ is the number of occurrences of $c^{\pm
1}$ in $A$, called the $c$-{\em length} of the word $A$. Conjugating by an
element of  $\PP_{B, X}$ if necessary, we can also assume that $A$ starts
with $c^{\pm 1}$.

Consider reduced diagrams $\Delta(j)$, $j=1,2,3$, over $\HH_c$ such that
$$
\p \Delta(j) = \ell(j) u(j) r(j) d(j) ,
$$
where $ \ph(\ell(j)) = p_{j-1}$, \ $\ph(r(j))^{-1} = p_{j}$, \  $\ph(u(j)) =
\ph(d(j))^{-1}  = A$, and $j=1,2,3$. Since $|\ell(j)|_c = |r(j)|_c = 0$ and
$|A|_c$ is minimal, it follows that every $c$-bond in $\Delta(j)$ connects an
edge on $u(j)$ to an edge on $d(j)$, $j=1,2,3$. Let
$$
\E_1(j), \dots, \E_{|A|_c}(j)
$$
be all consecutive along $u(j)$ $c$-bonds and let
$$
\p \E_i(j) = \ell_i(j) u_i(j)  r_i(j) d_i(j)
$$
be the standard contour of $\E_i(j)$, where $u_i(j) = \p \E_i(j) \wedge
u(j)$, \ \  $d_i(j) = \p \E_i(j) \wedge d(j)$ are the edges of $u(j)$, $d(j)$
labelled by $c^{\pm 1}$, $c^{\mp 1}$, respectively, and
\begin{gather*}
u(j) = u_1(j) v_1(j)  u_2(j) v_2(j) \dots u_{|A|_c}(j) v_{|A|_c}(j) , \\
d(j)^{-1} = d_1(j)^{-1} e_1(j)^{-1}   d_2(j)^{-1} e_2(j)^{-1}
\dots d_{|A|_c}(j)^{-1}  e_{|A|_c}(j)^{-1} ,
\end{gather*}
where $v_i(j)$ (resp. $e_i(j)^{-1}$)  is a subpath of $u(j)$ (resp.
$d(j)^{-1}$) which sits between sections $u_i(j)$ and $u_{i+1}(j)$  of
$u(j)$ (resp. sections $d_i(j)^{-1}$ and $d_{i+1}(j)^{-1}$  of  $d(j)^{-1}$),
$i=1, \dots, |A|_c-1$, $j=1,2,3$.

By $\Delta_i(j)$, $i=1, \dots, |A|_c$,   denote  the subdiagram of
$\Delta(j)$ which sits between $c$-bonds $\E_i(j)$ and $\E_{i+1}(j)$, where
$\E_{i+1}(j) = r(j)$ if $i = |A|_c$. Clearly, $|\p \Delta_i(j) |_c = 0$,
$\Delta_i(j)$ contains no $c$-edges and
$$
\p \Delta_i(j) = r_i(j)^{-1} v_i(j) \ell_{i+1}(j)^{-1} e_i(j) ,
$$
where $i = 1, \dots, |A|_c$ and $\ell_{i+1}(j)^{-1}  = r(j)$ if $i = |A|_c$.

Suppose that for some $i \in \{ 1, \dots, |A|_c -1 \}$ it is true that
\begin{equation}
\label{ui=ui1}
 \ph(u_i(1)) = \ph(u_{i+1}(1)) .
\end{equation}
Then, considering the diagram $\Delta_i(1)$, we can see that elements
$\ph(r_i(1))^{-1}$ and $\ph(\ell_{i+1}(1))$ are conjugate in the group
$\PP_{B, X}$ by $\ph(v_i(j))^{-1}$. It follows from definitions and equality
(\ref{ui=ui1}) that one of $r_i(1)^{-1}$, $ \ell_{i+1}(1)$ is a word in
$\X^{\pm 1} = \{ x_1^{\pm 1}, \dots,  x_m^{\pm 1} \}$ and the other one is a
word in the alphabet $\{ (b_1 x_1)^{\pm 1}, \dots,  (b_m x_m)^{\pm 1} \}$.
However, the conjugacy of such words in the direct product $\PP_{B, X} =
B(\B, n) \times B(\X, n)$ implies that $\ph(r_i(1)) =1$ whence $p_0 =1$. This
contradiction to $p_0 \neq 1$ enables us to conclude that $\ph(u_i(1)) =
\ph(u_{i+1}(1))^{-1}$ for every $i = 1, \dots, |A|_c-1$.

Switching to the diagram $\Delta_{|A|_c}(1)$ from  $\Delta_i(1)$, by a
similar argument we can conclude that $\ph(u_{|A|_c}(1))= \ph(u_1(1))^{-1}$.
In particular, this implies that  $|A|_c > 1$ and there is an $i \in \{1,
\dots, |A|_c \}$ such that
$$
\ph(u_i(1)) = c
$$
and either $i <  |A|_c$ or, otherwise, $i = |A|_c =2$.

Let $ \ph(\ell_i(1)) = W_1(\X)$, where $W_1(\X)$ is an $\X$-word. Then
$$
\ph(r_i(1))^{-1} = W_1(x_1 b_1, \dots,x_m b_m) = W_1(\B) W_1(\X) .
$$
Also, we let $ \ph(v_i(1)) = U(\B) V(\X)$, where $U(\B)$,  $V(\X)$ are $\B$-,
$\X$-words, respectively.

It is convenient to denote $r(1)$ by $\ell_{|A|_c+1}(1)^{-1}$, that is,
to define $\ell_{i+1}(1)$ for $i = |A|_c$.
Then $\ph( \ell_{|A|_c+1}(1) ) = \ph( \ell_1(2))$.

Considering the subdiagram $\Delta_i(1)$, we can see that
\begin{multline*}
\ph(\ell_{i+1}(1)) = \ph(v_{i}(1))^{-1}  \ph(r_i(1))^{-1}  \ph(v_{i}(1))  = \\
= (U(\B) V(\X))^{-1}  W_1(\B) W_1(\X) U(\B) V(\X)
\end{multline*}
in the group $\PP_{B, X}$. Since $\ph(\ell_{i+1}(1)) =  W_2(\B) W_2(\X)$
in $\PP_{B, X} = B(\B, n) \times B(\X, n)$  and both $B(\B, n)$, $B(\X, n)$
are free Burnside groups, it follows that
\begin{equation}
\label{W=W}
U(\X)^{-1} W_1(\X) U(\X) = V(\X)^{-1} W_1(\X) V(\X)
\end{equation}
in the group $B(\X, n)$. It follows from $p_0 \neq 1$ that $W_1(\X) \neq 1$ in
$B(\X, n)$.

Recall that if $n$ is odd and $n \gg 1$ (e.g., $n > 10^{10}$) then every
nontrivial element of $B(\X, n)$ is contained in a unique maximal cyclic
subgroup of $B(\X, n)$ of order $n$ (see Theorem 19.4 \cite{O89}) and every
maximal cyclic subgroup $C$ of $B(\X, n)$ is {\em antinormal} in $B(\X, n)$,
that is, for every $g \in B(\X, n)$ the inequality $g  C g^{-1} \cap C \neq
\{ 1\}$ implies that $g \in C$ (this is actually shown in process of the
proof of Theorem 19.4 \cite{O89}; similar arguments can also be found in
\cite{AI87}, \cite{I92a}).

Let $W_0(\X)$ be an $\X$-word so that $\langle W_0(\X) \rangle \subseteq
B(\X, n)$ is the maximal cyclic subgroup of order $n$ that contains
$W_1(\X)$. Since the subgroup $\langle W_0(\X) \rangle$ is antinormal in
$B(\X, n)$, we have from (\ref{W=W}) that $U(\X) V(\X)^{-1} \in \langle
W_0(\X) \rangle$.

Now consider the word $\ph( \ell_i(2) )$. It is easy to see that
$$
\ph( \ell_i(2)) = S(\X)^{-1} W_1(\X) S(\X) ,
$$
where $S(\X)$ is the product
$$
\ph( v_i(1)) \dots \ph( v_{|A|_c}(1))
\ph( v_1(2)) \dots \ph( v_{i-1}(2)) \in \PP_{B, X} .
$$

Replacing the word $\ph(\ell_i(1)) = W_1(\X)$ by
$\ph(\ell_i(2)) = S(\X)^{-1} W_1(\X) S(\X)$ and the diagram $\Delta(1)$ by
$\Delta(2)$  in the above argument, we will analogously obtain that
$$
U(\X)V(\X)^{-1}   \in  S(\X)^{-1} \langle W_0(\X) \rangle S(\X) .
$$

Since the maximal cyclic subgroups  $S(\X)^{-1} \langle W_0(\X) \rangle
S(\X)$ and $\langle W_0(\X) \rangle$ are either coincide or have the trivial
intersection, it follows that either $U(\X)V(\X)^{-1}=1$ in $B(\X, n)$ or
$S(\X) \in \langle W_0(\X) \rangle$. In the case $U(\X)V(\X)^{-1}=1$, we
distinguish two subcases: $i < |A|_c$ and  $i = |A|_c = 2$. If  $i < |A|_c$
then we have a contradiction to the minimality of $|A|_c$ because
$$
\ph( u_i(1)) \ph( v_i(1)) \ph(u_{i+1}(1)) = c U(\B) U(\X) c^{-1} = U(\X)
$$
in $\HH_c$ and $|U(\X)|_c = 0$. If  $i = |A|_c = 2$ then the word $A =
c^{-1}\ph( v_1(1)) c \ph( v_2(1))$ is conjugate in the group $\HH_c$ given by
(11) to an element of $\PP_{B, X}$ because
$$
c \ph( v_2(1)) = c U(\B) U(\X) = U(\X) c
$$
and so $A =  c^{-1} \ph( v_1(1)) U(\X) c$ in $\G_U$, contrary to our
assumption. These contradictions prove that the case  $U(\X)V(\X)^{-1}=1$ is
impossible and so $S(\X) \in \langle W_0(\X) \rangle$. Therefore,
$$
\ph(\ell_i(2)) = S(\X)^{-1} W_1(\X) S(\X) =  W_1(\X)
$$
in $\G_U$, whence $A^{-1}  p_0 A = p_0$. Now we can conclude that
$p_0 \in \F(A)$, $\kp_A = 1$ in $\Aut \F(A)$ and the group
$\langle \kp_A, A \rangle$ has   exponent $n$, as required.
Lemma 2 is proven. \qed
\smallskip

By results of \cite{I02b} and Lemma 2, the group $\PP_{B, X} = B(\B, n)
\times B(\X, n) $ naturally embeds in the quotient $\HH_c / \HH_c^n$. Denote
the subgroup $\langle \B, c \rangle$ of $\HH_c / \HH_c^n$, generated by $\B$
and $c$, by
\begin{equation}
\label{KBc}
\K_{B,c} = \la B, c \ra \subset \HH_c / \HH_c^n .
\end{equation}
Observe that, in view of relations $x_i b_j  x_i^{-1}  = b_j$, \  $x_i c
x_i^{-1} = c b_i$,  \ $i, j \in \{ 1, \dots, m \}$, which hold in $\HH_c /
\HH_c^n$,  the group $\HH_c / \HH_c^n$ is a splitting extension of $\K_{B,
c}$ by $B(\X, n)$. In particular, we can also consider the following
HNN-extension of $\K_{B, c}$
\begin{multline}
\label{HHX} \HH_X = \langle \K_{B, c}, x_1, \dots, x_m  \; \| \; x_i b_j
x_i^{-1}  = b_j,  x_i c x_i^{-1} = c b_i, \; i,j \in \{ 1, \dots, m \}
\rangle .
\end{multline}
and claim that the quotient $\HH_X / \HH_X^n$ is naturally isomorphic to
$\HH_c / \HH_c^n$.

\begin{lemma} Every maximal cyclic subgroup of the group
$\K_{B, c}$ is antinormal in $\K_{B, c}$.
\end{lemma}

{\em Proof.} Consider a graded presentation for the group $\HH_c / \HH_c^n$
which, by results of \cite{I02b} and Lemma 2, can be constructed as in
\cite{I02b}:
\begin{multline}
\label{HHcn}
\HH_c / \HH_c^n =  \HH_c(\infty)  = \langle \PP_{B, X}, c  \; \| \; c^{-1}
x_j c =  x_j b_j, \; j \in \{ 1, \dots, m\},
\\
A^n=1, A \in \cup_{i=1}^\infty \X_i \rangle .
\end{multline}

Let $A$ be a period of rank $i \ge 1$, that is, $A \in \X_i$. It is easy to
see from the definition of $\X_i$ and details of the proof of Lemma 2 that
$\kp_A = 1$ in $\Aut \F(A)$ (that is, $A$ always centralizes $\F(A)$) and
that the group $\F(A)$ consists of words either of the form $F(x_1, \dots,
x_m)$ (if $A$ starts with $c$; recall that, by definitions of \cite{I02b},
$A$ starts with $c^{\pm 1}$) or of the form $F(x_1 b_1, \dots, x_m b_m)$ (if
$A$ starts with $c^{-1}$).

Since there is a natural homomorphism $\HH_c(\infty) \to B(\X, n)$ whose kernel is
$\K_{B, c}$, it follows from the above observation  that if $F \in \F(A)$ or $p \in
\PP_{B, X}$ is conjugate in $\HH_c(\infty)$ to an element of $\K_{B, c}$ then $F =1$ or
$p \in B(\B, n)$, respectively. We will repeatedly use this remark in the arguments
below.

Suppose that $\Ss$ is a subgroup of $\langle A, \F(A) \rangle \subseteq
\HH_c(\infty) $, where $A \in \X_i$, and $\Ss$ is conjugate in
$\HH_c(\infty)$ to a subgroup of $\K_{B, c}$. By the above remark, every
nontrivial element of  $\Ss$ has the form $A^k F$, where $k \not\equiv 0
\pmod n$ and $F \in \F(A)$. In particular, if $A^kF \in \Ss$ and $kd \equiv 0
\pmod n$ then
$$
(A^k F)^d = A^{k d} F^d = F^d = 1 .
$$
Similarly, if $A^{k_1}F_1,  A^{k_2}F_2 \in \Ss$  and $k_1 \equiv k_2 \pmod n$
then $F_1 = F_2$ in $\F(A)$. Therefore, we can conclude that $\Ss$ is cyclic
(and is generated by an element of the maximal order).

Now assume that $\langle W_0 \rangle$ is a maximal cyclic subgroup of $\K_{B, c}$, $W_1
\in \K_{B, c}$ and
\begin{equation}
\label{W0}
W_1 W_0^{\ell_1} W_1^{-1} = W_0^{\ell_2}
\end{equation}
in $\HH_c(\infty)$, where $W_0^{\ell_2} \neq 1$. By Lemma 10.2 \cite{I02b},
$W_0$ is conjugate in $\HH_c(\infty)$ to either an element $U_0$ with
$|U_0|_c = 0$ (by the above remark, $U_0 \in B(\B, n)$) or to $A^{k_0}F_0$,
where $A$ is a period of rank $i \ge 1$, $0 < k_0 < n$, and $F_0 \in \F(A)$.

Conjugating, we may assume that $W_0 =U_0$ or $W_0 = A^{k_0}F_0$ as above
(and $\la W_0, W_1 \ra$ is conjugate in $\HH_c(\infty)$ to a subgroup of
$\K_{B, c}$). Consider a disk diagram $\Delta$ for equality (\ref{W0}) and
let  $\Delta_0$  be an annular diagram obtained from $\Delta$ by identifying
sections of the contour of $\Delta$ labelled by $W_1$ and $W_1^{-1}$. Since
contours of $\Delta_0$  are labelled by $W_0^{\ell_1} = U_0^{\ell_1}$ or
$W_0^{\ell_1} = A^{k_0 \ell_1} F_0^{\ell_1}$ and $W_0^{-\ell_2} =
U_0^{-\ell_2}$ or $W_0^{-\ell_2} = A^{-k_0 \ell_2} F_0^{-\ell_2}$, it follows
that every $c$-annulus in $\Delta_0$  will be contractible. Otherwise, as in
the above remark, a word of the form $F(x_1, \dots, x_m)$ or $F(x_1 b_1,
\dots, x_m b_m)$ would be conjugate in $\HH_c(\infty)$ to $W_0^{\ell_1}$,
whence $W_0^{\ell_1}= 1$ which is a contradiction. Hence, we can bring
$\Delta_0$  to a reduced form $\Delta_1$ by removal of reducible pairs  and
reducible $c$-annuli. (Recall that $n$ is odd and so there are no
self-compatible 2-cells, see \cite{I02b}.)

Hence, there is a path in $\Delta_1$ (like in $\Delta_0$)
which connects vertices on distinct contours
of $\Delta_1$ and has label equal in $\HH_c(\infty)$ to $W_1$.
Now we can argue as  in the proof of Lemma 10.2 \cite{I94}
to show that if $W_0 = A^{k_0}F_0$ then $W_1 = A^{k_1}F_1$  in $\HH_c(\infty)$,
where $F_1 \in \F(A)$, and if $W_0 = U_0$ then
$W_1 = U_1$  in $\HH_c(\infty)$, where $U_1 \in \PP_{B,X}$,
and $U_1 U^{\ell_1}_0 U_1^{-1} = U^{\ell_2}_0$ in $B(\B, n)$.

In the latter case $W_0 = U_0$,  we note that $U_1 \in B(\B, n)$ by the above
remark and, by the antinormality of maximal cyclic subgroups of $B(\B, n)$
(see the proof of Lemma 2), we have that the subgroup $\la U_0, U_1 \ra$  is
cyclic. Now the required inclusion $U_1 \in \la U_0 \ra$ follows from the
original maximality of $\la U_0 \ra$ in $\K_{B, c}$.

In the former case $W_0 = A^{k_0}F_0$, the subgroup $ \langle  A^{k_0}F_0,
A^{k_1}F_1 \rangle \subseteq \langle  A, \F(A)  \rangle $ is conjugate to a
subgroup of $\K_{B, c}$ and, by the above remark, is cyclic. As before, it
follows from the original maximality of $\langle W_0 \rangle$ in $\K_{B, c}$
that $A^{k_0}F_0$ is a generator of $\langle A^{k_0}F_0, A^{k_1}F_1  \rangle
$. Hence, we have again that $W_1 \in \langle W_0 \rangle$,  as required. The
proof of Lemma 3 is complete. \qed
\smallskip

Replacing the letters $b_i$, $x_i$, $c$ in the above construction by $\bar
b_i$, $\bar x_i$, $\bar  c$, $i=1, \dots, m$, we analogously to (\ref{HHc}),
(\ref{KBc}), (\ref{HHX}) construct groups $\HH_{\bar c}$, $\K_{\bar B, \bar
c}$, $\HH_{\bar X}$, respectively.

Consider an HNN-extension $\G_{X, \bar X}$ of the direct product $\K_{B, c}
\times \K_{\bar B, \bar c}$ with stable letters $x_1, \bar x_1, \dots$, $x_m,
\bar x_m$ defined as follows.
\begin{multline}
\label{HHXX} \HH_{X, \bar X} = \langle \K_{B, c} \times \K_{\bar B, \bar c}, x_1, \bar
x_1, \dots, x_m, \bar x_m \; \| \;
x_i b_j  x_i^{-1}  = b_j,  \bar x_i \bar b_j  \bar x_i^{-1}  = \bar b_j, \\
x_i c x_i^{-1} = c b_i, \bar x_i \bar c \bar x_i^{-1} = \bar c \bar b_i, \; i,j \in \{ 1,
\dots, m \} \rangle .
\end{multline}

\begin{lemma} Property $(E)$ holds for the group $\HH_{X, \bar X}$,
given by $(\ref{HHXX})$, with $N_E = 1$.
\end{lemma}

{\em Proof.} Let $A \in \HH_{X, \bar X}$, $A$ not conjugate to an element $p
\in \K_{B, c} \times \K_{\bar B, \bar c}$, and $A^{-1} p_0 A = p_1$, where
$p_0, p_1 \in \K_{B, c} \times \K_{\bar B, \bar c}$. We may assume that if $A
= B$ in $\G_{X, \bar X}$ then $|A|_{\X \cup \bar \X} \le |B|_{\X \cup \bar
\X}$ ($|A|_{\X \cup \bar \X}$ is the ${\X \cup \bar \X}$-length of $A$, that
is, the number of all occurrences of $x^{\pm 1}$ where $x \in \X \cup \bar
\X$) and that $A$ starts with $x^{\pm 1}$, where $x \in \X \cup \bar \X$.

Considering a diagram $\Delta$ over the group   $\HH_{X, \bar X}$ with $\p
\Delta = \ell u r d$, where
$$
\ph(\ell) = p_0, \  \ph(u) = A, \
\ph(r) = p_1^{-1}, \  \ph(d) = A^{-1} ,
$$
and arguing as in the proof of Lemma 3, we can show that all $\X \cup \bar
\X$-letters of $A$ are either in $\X^{\pm 1}$ or in $\bar \X^{\pm 1}$
(otherwise, $p_0 = 1$). For definiteness, assume that all $\X \cup \bar
\X$-letters of $A$ are in $\X^{\pm 1}$ (the other case is symmetric). Then it
is clear that $\F(A) = \K_{B, c}$ and $p_0 \in \K_{B, c}$. Let $\wtl A$
denote the image of $A$ under the projection $\HH_{X, \bar X} \to \HH_{X}$,
where $\HH_{X}$ is given by presentation  (\ref{HHX}), which erases all
letters $\bar x_j$, $\bar b_j$, $\bar c$, $j = 1,\dots, m$. Note that $A$
acts on $\F(A)$ by conjugations as the word $\wtl A$ does. Therefore, the
group $\la \kp_A, \F(A) \ra$ is a homomorphic image of the subgroup $\la \wtl
A, \F(A) \ra \subseteq \HH_{X}/ \HH_{X}^n$ (recall that $\F(A) = \K_{B, c}$
and $\K_{B, c}$ embeds in $\HH_{X}/ \HH_{X}^n$).  Thus, the group $\la \kp_A,
\F(A) \ra$ has exponent $n$ and so property (E) holds for $A$ with $N_E =1$.
Lemma 4 is proved. \qed
\smallskip

Making the substitution $ d_i \to b_i$,  \ $y_i \to x_i$, \  $e \to c$, \ $\bar  d_i \to
\bar b_i$, \ $\bar y_i \to \bar x_i$,  \ $\bar e \to \bar c$, we construct groups
$$
\HH_e , \ \ \HH_{\bar e} , \ \  \K_{D, e} ,  \ \ \K_{\bar D, \bar e} , \ \ \HH_{Y, \bar
Y}
$$ which are complete $(d, y, e)$-analogs of corresponding $(b, x, c)$-groups. In
particular, we define
\begin{multline}
\label{HHYY}
\HH_{Y, \bar Y} = \langle \K_{D, e} \times \K_{\bar D, \bar e}, y_1, \bar
y_1, \dots, y_m, \bar y_m \; \| \; y_i d_j  y_i^{-1}  = d_j,  \bar y_i \bar
d_j  \bar y_i^{-1}  = \bar d_j, \\ y_i e y_i^{-1} = e b_i, \bar y_i \bar e
\bar y_i^{-1} = \bar e \bar d_i, \; i,j \in \{ 1, \dots, m \} \rangle .
\end{multline}

Observe that the subgroup $\la b_1 \bar b_1, \dots, b_m \bar b_m, c \bar c
\ra$ of $\K_{B, c} \times \K_{\bar B, \bar c}$ is naturally isomorphic to
$\K_{B, c}$ and the subgroup $\la d_1 \bar d_1 a_1, \dots, d_m \bar d_m a_m,
e \bar e \ra$ of $\K_{D, e} \times \K_{\bar D, \bar e}  \times \G_A/\G_A^n$
is also naturally isomorphic to $\K_{D, e}$ (which is naturally isomorphic to
$\K_{B, c}$), because the subgroup $\la d_1 \bar d_1, \dots, d_m \bar d_m
\ra$ of  $\K_{D, e} \times \K_{\bar D, \bar e}$ is naturally isomorphic to
the free Burnside group $B(\D, n)$ and so for every $\A$-word $W$ an
inequality $W(a_1, \dots, a_m) \neq 1$ in $\G_A/\G_A^n$ implies that $W(  d_1
\bar d_1, \dots, d_m \bar d_m) \neq 1$ in $\K_{D, e} \times \K_{\bar D, \bar
e}$. Therefore, the map
$$
b_1\bar b_1 \to d_1\bar d_1 a_1 , \  \dots, \ b_m \bar b_m \to d_m \bar d_m a_m, \ \ c
\bar c \to d \bar d
$$
extends to an isomorphism of corresponding subgroups of $\K_{B, c} \times
\K_{\bar B, \bar c}$ and $\K_{D, e} \times \K_{\bar D, \bar e} \times
\G_A/\G^n_A$, where the group $\G_A$ is defined by (\ref{GA}).

Consider the free Burnside $n$-product
\begin{equation}
\label{Q} \Q = ( \HH_{X, \bar X} /  \HH_{X, \bar X}^n )  *_n ( \HH_{Y, \bar
Y} / \HH_{Y, \bar Y}^n  \times \G_A / \G_A^n )
\end{equation}
of the groups $\HH_{X, \bar X} /  \HH_{X, \bar X}^n $ and $\HH_{Y, \bar Y} /
\HH_{Y, \bar Y}^n  \times \G_A / \G_A^n $ (which is the quotient $F_P /F_P^n$
of their free product $F_P = (\HH_{X, \bar X} /  \HH_{X, \bar X}^n )  * (
\HH_{Y, \bar Y} /  \HH_{Y, \bar Y}^n  \times \G_A / \G_A^n )$ by $F_P^n$).
Next, in view of the isomorphism pointed out above, we can define an
HNN-extension $\HH_z$ of $\Q$  as follows:
\begin{equation}
\label{HHz} \HH_z = \langle \Q, z  \; \| \; z b_i \bar b_i z^{-1} = d_i \bar
d_i a_i, \; z c \bar c z^{-1} = e \bar e, \; i \in \{ 1, \dots, m \} \rangle
.
\end{equation}

\begin{lemma} Property $(E)$ holds for the group $\HH_{z}$,  given
by $(\ref{HHz})$, with $N_E = 3$.
\end{lemma}

{\em Proof.} This proof is analogous to the proof of Lemma 2.
Assume that $A \in \HH_z$, $A$ is not conjugate in $\HH_z$ to $p \in \Q$
and $A^{-1}  p_{j-1} A = p_j$, where $p_{j-1}, p_j \in \Q$ and $j=1,2,3$,
$p_0 \neq 1$. Our goal is to show that $p_0 \in \F(A)$ and the subgroup
$\langle \kp_A, \F(A) \rangle$ has exponent $n$.

As in the proof of Lemma 2, we can assume that $A$ is a reduced in $\HH_z$ word,
that is,  if $A = B$ in $\HH_z$ then
$|A|_z \leq |B|_z$, where $|A|_z$ is the $z$-{\em length} of the word $A$.
Conjugating if necessary, we can also
assume that $A$ starts with $z^{\pm 1}$.

As in the proof of Lemma 2, we consider reduced diagrams $\Delta(j)$,
$j=1,2,3$, over $\HH_z$ and introduce the notation related to these diagrams
and their $z$-bonds exactly as in that proof. As before, if for some $i <
|A|_z$ we had that $ \ph(u_i(1)) = \ph(u_{i+1}(1)) $ then, by the definition
of presentation (\ref{HHz}) of $\HH_z$, we would have that one of words
$\ph(r_i(1))^{-1}$, $\ph(\ell_{i+1}(1))$  would be a word in the alphabet $\{
b_1 \bar b_1, \dots, b_m \bar b_m, c \bar c \}^{\pm 1}$ and the other one
would be in the alphabet $\{ d_1 \bar d_1 a_1, \dots, d_m \bar d_m a_m, e
\bar e \}^{\pm 1}$. However, these words are conjugate in $\Q$ and so they
must represent the trivial element of $\Q$, contrary to $p_0 \neq 1$. By a
similar argument, we can conclude that $\ph(u_{|A|_z}(1))= \ph(u_1(1))^{-1}$.
In particular, as in the proof of Lemma 2, this implies that $|A|_z > 1$ and
there is an $i \in \{ 1, \dots , |A|_z \}$  such that
$$
\ph(u_i(1)) = z
$$
and either $i < |A|_z$ or, otherwise, $i = |A|_z =2$. Let
$$
\ph(r_i(1))^{-1} = W_1(  b_1 \bar b_1, \dots, b_m \bar b_m, c \bar c  ) .
$$
As in before, it is convenient to denote $r(1)$ by $\ell_{|A|_c+1}(1))^{-1}$,
that is, to define $\ell_{i+1}(1)$ for $i = |A|_z$.
Then $\ph( \ell_{|A|_z+1}(1) ) = \ph( \ell_1(2))$.

Recall that factors $F_1, F_2$ of a free Burnside $n$-product $F_P / F_P^n =
F_1 *_n F_2$ of odd exponent $n \gg 1$ (say, $n > 10^{10}$), where $F_P = F_1
* F_2$ is the free product, are antinormal in $F_P/ F_P^n$ (e.g., see  Lemma
34.10 \cite{O89}). Since
\begin{equation}
\label{22}
\ph(v_i(1)) \ph(\ell_{i+1}(1)) \ph(v_i(1))^{-1} = \ph(r_i(1))^{-1} \neq 1
\end{equation}
and  $\ph(\ell_{i+1}(1)) $, $\ph(r_i(1))^{-1} $ are in $\HH_{X, \bar X}/
\HH_{X, \bar X}^n$, it follows that $\ph(v_i(1))$ is also in  $\HH_{X, \bar
X}/ \HH_{X, \bar X}^n$. Hence, we may suppose that  $\ph(v_i(1))$ is a word
in the alphabet
$$
\B^{\pm 1} \cup \bar \B^{\pm 1} \cup \{ c^{\pm 1}, \bar c^{\pm 1} \} \cup \X^{\pm 1} \cup
\bar \X^{\pm 1} .
$$

Consider a graded presentation for the group $\HH_{X, \bar X}/ \HH_{X, \bar
X}^n$ which, by results of \cite{I02b} and Lemma 4, can be constructed as in
\cite{I02b}:
\begin{multline}
\label{HHXXn} \HH_{X, \bar X} / \HH_{X, \bar X}^n  =
\HH_{X, \bar X}(\infty)= \langle \;
\K_{B, c} \times \K_{\bar B, \bar c}, x_1, \bar x_1,
\dots, x_m, \bar x_m \; \| \; x_j b_k  x_j^{-1}  = b_k, \\
\bar x_j \bar b_k  \bar x_j^{-1}  = \bar b_k, \; x_j c x_j^{-1} = c b_j, \;
\bar x_j \bar
c \bar x_j^{-1} = \bar c \bar b_j, \; j,k \in \{ 1, \dots, m \},   \\
A^n =1, A \in \cup_{i=1}^\infty \X_i \;  \rangle .
\end{multline}

Let $\Delta$ be a disk diagram over the presentation (\ref{HHXXn}) for the
equality (\ref{22}) and $\Delta_0$ be an annular  diagram obtained from
$\Delta$ by identifying sections of the contour of $\Delta$ labelled by
$\ph(v_i(1))$ and $\ph(v_i(1))^{-1} $. Since contours of $\Delta_0$ are
labelled by $W_1(  b_1 \bar b_1, \dots, b_m \bar b_m, c \bar c  )$ and
$$
\ph(\ell_{i+1}(1))^{-1}= W_2(  b_1 \bar b_1, \dots, b_m \bar b_m, c \bar c  ) ,
$$
it follows that every $x$- or $\bar x$-annulus will be contractible
(otherwise,  a word of the form  $F(b_1, \dots, b_m, c)$ or $F(\bar b_1,
\dots, \bar b_m,\bar c )$ would be conjugate in $\K_{B, c} \times \K_{\bar B,
\bar c}$ to $W_1(  b_1 \bar b_1, \dots, b_m \bar b_m, c \bar c)$ or to $W_2(
b_1 \bar b_1, \dots, b_m \bar b_m, c \bar c)^{-1}$, whence $p_0 = 1$). This
observation (together with a standard argument as in the proof of Lemma 3)
enables us to conclude that if $\ph(u_i(1)) = z$ then
\begin{equation}
\label{inB} \ph(v_i(1))  \in \K_{B, c} \times \K_{\bar B, \bar c} .
\end{equation}

Now suppose that
$$
\ph(u_i(1)) = z^{-1}
$$
for some $i \in \{ 1, \dots , |A|_z \}$. Then $\ph(r_i(1))^{-1}$,
$\ell_{i+1}(1)$ are words in the alphabet $\{ d_1 \bar d_1 a_1, \dots, d_m
\bar d_m a_m, e \bar e \}^{\pm 1}$ and, in the above arguments, we have to
replace the group $\HH_{X, \bar X} / \HH_{X, \bar X}^n$ by $ (\HH_{Y, \bar
Y}/  \HH_{Y, \bar Y}^n ) \times  ( \G_A / \G_A^n )$. However, if we disregard
$\G_A$-components of elements of the direct product $ (\HH_{Y, \bar Y}/
\HH_{Y, \bar Y}^n ) \times  ( \G_A / \G_A^n )$ then our arguments are
retained and so we can  conclude that if $\ph(u_i(1)) = z^{-1}$ then
\begin{equation}
\label{inD} \ph(v_i(1)) \in   \K_{D, e} \times \K_{\bar D, \bar e}
\times  \G_A / \G_A^n
.
\end{equation}

Now suppose again that $\ph(u_i(1)) = z$ for some $i \in \{ 1, \dots , |A|_z
\}$ (and either $i < |A|_z$ or $i = |A|_z= 2$). Denote
\begin{gather*}
\ph(r_{i}(1))^{-1}=
W_1(  b_1 \bar b_1, \dots, b_m \bar b_m, c \bar c  ) =W_1 \bar W_1 , \\
\ph(\ell_{i+1}(1))=
W_2( b_1 \bar b_1, \dots, b_m \bar b_m, c \bar c  ) = W_2 \bar W_2 , \\
\ph(v_{i}(1)) = U(b_1, \dots, b_m, c ) V(\bar b_1, \dots, \bar b_m, \bar c  )
= U \bar V ,
\end{gather*}
where $W_1 = W_1(b_1 , \dots, b_m , c)$, $\bar W_1 =
W_1(\bar b_1, \dots, \bar b_m, \bar
c)$ and so on. Clearly,
$$
U^{-1} W_1 U = W_2 , \quad  \bar  V^{-1} \bar  W_1 \bar V =
\bar W_2 , \quad V^{-1}  W_1
V = W_2 ,
$$
whence
\begin{equation}
\label{VUW}
 VU^{-1} W_1 U V^{-1} = W_1
\end{equation}
in $\K_{B, c}$. Let the word
$$
P(\B, c, \bar \B, \bar c) =  P_1(b_1 , \dots, b_m , c) P_2( \bar b_1 , \dots,
\bar b_m , \bar c) = P_1 \bar P_2 \in \K_{B, c} \times  \K_{\bar B, \bar c}
$$
be obtained from the cyclic permutation $\bar A$ of $A = \ph(u(1))$ that
starts with $\ph(u_i(1))$ by deletion of all letters $z^{\pm 1}$, $a_j^{\pm
1} \}$ and by replacing the letters $d_j$, $\bar d_j$, $e$, $\bar e$ by
$b_j$, $\bar b_j$, $c$, $\bar c$, respectively, where $j =1, \dots, m$. It is
not difficult to see from inclusions (\ref{inB})--(\ref{inD}) that
$\ph(r_{i}(2))^{-1} = P^{-1} \ph(r_{i}(1))^{-1} P$ in $\K_{B, c} \times
\K_{\bar B, \bar c}$. Hence, we have in the group $\K_{B, c} \times \K_{\bar
B, \bar c}$ that
\begin{multline*}
\ph(r_{i}(2))^{-1} = P^{-1} \ph(r_{i}(1))^{-1} P =  P_1^{-1} W_1 P_1 \bar
P_2^{-1} \bar W_1 \bar P_2 = P_1^{-1} W_1 P_1   \bar P_1^{-1} \bar W_1 \bar
P_1 .
\end{multline*}

Repeating the above argument for $\ph(r_{i}(2))$, $\Delta(2)$, we analogously
to (\ref{VUW}) obtain that
\begin{equation}
\label{VUW2}
VU^{-1} P_1^{-1} W_1 P_1 U V^{-1} = P_1^{-1} W_1  P_1
\end{equation}
in $\K_{B, c}$. Let $\langle W_0 \rangle$ be a maximal cyclic subgroup of
$\K_{B, c}$ that contains $W_1 \neq 1$. By Lemma 3, it follows from
(\ref{VUW})--(\ref{VUW2}) that $V U^{-1} \in \langle W_0 \rangle$  and $V
U^{-1} \in P_1^{-1} \langle W_0 \rangle P_1$. By Lemma 3, this implies that
either $V U^{-1} =1$ in $\K_{B, c}$ or $\langle W_0 \rangle = P_1 \langle W_0
\rangle P_1^{-1}$. As in the proof of Lemma 2, the equality $V U^{-1} =1$
contradicts either the minimality of $|A|_z$ when $i < |A|_z$ or the
assumption that $A$ is not conjugate to an element $p \in \Q$ when $i = |A|_z
= 2$. Therefore, $\langle W_0 \rangle = P_1^{-1} \langle W_0 \rangle P_1$
and, by Lemma 3, $P_1 \in \langle W_0 \rangle$. Now we can conclude that
$A^{-1} p_0 A = p_0$, $p_0 \in \F(A)$, $\kp_A = 1$ and $\la \kp_A, \F(A) \ra$
has exponent $n$. The proof of Lemma 5 is complete. \qed

\section{Proofs of Theorems 1--3}

{\em Proof of Theorem 1.} Observe that, by definition of the presentation
(\ref{HHz}) for $\HH_z$, there is a natural homomorphism $\G_U \to \HH_z$.
Indeed, defining relations (\ref{R7}) of $\G_U$ hold in $\HH_z$ given by
(\ref{HHz}), relations (\ref{R6}) of $\G_U$ hold in $\Q$ given by (\ref{Q}),
relations (\ref{R1}), (\ref{R2}), (\ref{R4}), (\ref{R5}) of $\G_U$ hold in
groups $\HH_{X, \bar X}$, $\HH_{Y, \bar Y}$ given by (\ref{HHXX}),
(\ref{HHYY}) and relations (\ref{R3}) of $\G_U$ hold in groups $\K_{B, c}$,
$\K_{\bar B, \bar c}$, $\K_{D, e}$, $\K_{\bar D, \bar e}$ (see (\ref{KBc})).

Furthermore, note that defining relations of the group $\HH_z$ are
(\ref{R1})--(\ref{R7}) and all other defining relators of $\HH_z$ are words
in $\Ker \alpha_U \cdot F(\U)^n$, where $\alpha_U : F(\U) \to \G_U$ is the
natural homomorphism (because this is the case for groups  $\Q$, $\HH_{X,
\bar X}$, $\HH_{Y, \bar Y}$).  Consequently, we may conclude that the groups
$\G_U / \G_U^n$ and $\HH_z / \HH_z^n$ are naturally isomorphic. By Lemma 5
and results of \cite{I02b}, the group $\Q$ and hence $\G_A / \G_A^n = \la
a_1, \dots, a_m \ra \subseteq \Q$ both naturally embed in the quotient $
\HH_z / \HH_z^n = \G_U / \G_U^n$.

Let $n > 2^{48}$ be odd and
$$
\G_{A,0} = \la \A \; \| \; R=1, R \in \R_{A, 0} \ra ,
$$
where $\R_{A, 0} = F(\A)^n$. Then $\G_{A,0} = F(\A) / F(\A)^n = B(\A,n)$ is
the free Burnside group of exponent $n$ in $\A$.  By Lemma 1, the group
$\G_{A,0} = \la a_1, \dots, a_m \ra$ embeds in $\G_U = \G_U(\G_{A,0})$ given
by (\ref{GU}) and, by Lemma 5 and results of \cite{I02b}, $\G_{A,0} = \la
a_1, \dots, a_m \ra$ also naturally embeds in $\G_U(\G_{A,0}) /
\G_U(\G_{A,0})^n$.

Let $\K$ be a normal subgroup of  $\G_{A,0}$ and
$$
\G_{A,1} =  \G_{A,0} / \K = \la \A \; \| \; R=1, R \in \R_{A, 1} \ra
$$
be a presentation for the quotient $\G_{A,1} =  \G_{A,0} / \K$. Consider
another group $\G_U = \G_U(\G_{A,1})$ defined by presentation (\ref{GU}) in
which we now use the group $\G_{A,1}$.  By Lemma 1, results of \cite{I02b}
and Lemma 5, the group $\G_{A,1} = \la a_1, \dots, a_m \ra$  naturally embeds
in both $\G_U(\G_{A,1})$ and $\G_U(\G_{A,1}) / \G_U(\G_{A,1})^n$. On the
other hand, it is immediate from definitions that the group $\G_U(\G_{A,1})$
is naturally isomorphic to the quotient $\G_U(\G_{A,0}) / \la \K \ra^{
\G_U(\G_{A,0}) }$. Hence, $\la \K \ra^{ \G_U(\G_{A,0}) } \cap  \G_{A,0} = \K$
for an arbitrary normal subgroup $\K$ of $\G_{A,0}$ and so $\G_{A,0}$ is a
$Q$-subgroup of $\G_U(\G_{A,0})$.

Analogously, the group $\G_U(\G_{A,1}) / \G_U(\G_{A,1})^n$ is naturally
isomorphic to the quotient
$$
( \G_U(\G_{A,0}) /  \G_U(\G_{A,0})^n ) / \la \K \ra^{ \G_U(\G_{A,0}) /
\G_U(\G_{A,0})^n  }
$$
and, therefore, in view of the natural embedding $\G_{A,1} \to \G_U(\G_{A,1})
/ \G_U(\G_{A,1})^n$, we have that $\la \K \ra^{ \G_U(\G_{A,0}) /
\G_U(\G_{A,0})^n  } \cap \G_{A,0} = \K$ and so $\G_{A, 0} = B(\A,n)$ is a
$Q$-subgroup of $\G_U(\G_{A,0}) / \G_U(\G_{A,0})^n$ as well.

Thus, $\G_{A,0} = \la a_1, \dots, a_m \ra \to  \G_U(\G_{A,0}) $ is a required
embedding and Theorem 1 is proved. \qed
\smallskip

{\em Proof of Theorem 2.}  As in the proof of Theorem 1, let
$$
\G_{A,0} = \la \A \; \| \; R=1, R \in \R_{A, 0} \ra ,
$$
where  $\R_{A, 0} = F(\A)^n$, so that $\G_{A,0} = B(\A,n)$. Let a verbal
subgroup $\V(B(2,n))$ of $B(2,n)$ be nontrivial. By results of \cite{I02c},
we can pick $|\U|+2$ $\A$-words $V_1, \dots, V_{|\U|}$, $V_{|\U|+1}$, $V_{
|\U|+2 }$ in the verbal subgroup $\V( \G_{A,0}) \neq \{ 1\}$  which generate
a $Q$-subgroup of $\G_{A,0} = B(\A,n)$ naturally isomorphic to a free
Burnside group $B(|\U|+2, n)$ of rank $|\U|+2$ and exponent $n$. Note that
the property of being a $Q$-subgroup is obviously transitive. Hence, it
follows from  Theorem 1 that $\la V_1, \dots, V_{|\U|} \ra $ is a
$Q$-subgroup of both $\G_U(\G_{A,0})$ and $\G_U(\G_{A,0}) /
\G_U(\G_{A,0})^n$.

Let $(u_1,  \dots, u_{| \U |})$ be an ordering of all letters of the alphabet
$\U$ given by (\ref{U}) so that $u_1 = a_1, \dots$, $u_m = a_m$. By induction
on $i \ge 1$, define the group $\G_{A,i}$  by presentation
$$
\G_{A,i} = \la \A \; \| \; R=1, R \in \R_{A, i} \ra ,
$$
where $ \R_{A, i-1} \subseteq \R_{A, i} $ and the set difference $ \R_{A, i}
\setminus \R_{A, i-1}$ consists of all words of the form $W(V_1, \dots, V_{
|\U| } )$ such that $W(u_1, \dots, u_{ |\U| } ) =1$ in the group $\G_U(
\G_{A, i-1} )$ given by (\ref{GU}) with $\G_A =  \G_{A, i-1}$, $i \ge 1$.

It follows from definitions that   a group  $\G_U( \G_{A} )$, defined by
(\ref{GU}), can be given by defining relations (\ref{R1})--(\ref{R7}) and
defining relations (modulo $F(\A)^n$)  of $\G_A$. Therefore, without altering
the groups $\G_U(\G_{A, i})$, we can change sets $\R_{A, i}$ by $\R'_{A, i}$
so that $\R'_{A, 0}$ is empty,  $\R'_{A, 1}$ contains only relators
corresponding to relations (\ref{R1})--(\ref{R7}), and $ \R'_{A, i} \setminus
\R'_{A, i-1}$, when $i
>1$, consists of words of the form $W(V_1, \dots, V_{m } )$ such that
$W(a_1, \dots, a_{ m } ) =1$ in $\G_U( \G_{A, i-1} )$ or in $\G_{A, i-1}$
(for $\G_{A, i-1}$ naturally embeds in $\G_U( \G_{A, i-1} )$ by Lemma 1).

Consider the limit group
$$
\G_{A,\infty} = \la \A \; \| \; R=1, R \in \cup_{i=0}^\infty \R'_{A, i} \ra  .
$$
It follows from definitions that if $W(u_1, \dots, u_{ |\U| } ) = 1$ in the
group $\G_U(\G_{A,\infty})$, defined by (\ref{GU}) with $\G_{A} =
\G_{A,\infty}$,  then $W(V_1, \dots, V_{ |\U| } ) = 1$ in $\G_{A,\infty}$.
Hence, the map $u_j \to V_j$, $j=1, \dots, |\U|$, extends to a homomorphism
$\G_U( \G_{A,\infty} ) \to \G_{A,\infty}$. By Theorem 1, the group
$\G_{A,\infty} $ naturally embeds in $ \G_U( \G_{A,\infty} ) / \G_U(
\G_{A,\infty} )^n$. On the other hand, it is easy to see from the definition
of defining relations of $\G_{A,\infty}$ that the map $u_j \to V_j$, $j=1,
\dots, |\U|$, extends to  an isomorphism
$$
\G_U( \G_{A,\infty} ) / \G_U( \G_{A,\infty} )^n \to \la V_1, \dots, V_{ |\U|
}  \ra \subseteq  \G_{A,\infty} .
$$
Hence, we can consider the following HNN-extension  of $ \G_U( \G_{A,\infty}
) /  \G_U( \G_{A,\infty} )^n$:
\begin{equation}
\label{HHt}
\HH_t = \la \G_U( \G_{A,\infty} ) /  \G_U( \G_{A,\infty} )^n, \; t \; \| \;
t u_j t^{-1} = V_j, \; j=1, \dots, |\U| \ra .
\end{equation}
Recall that  the set $\cup_{i=0}^\infty \R'_{A, i}$ consists of words in $\{
V_1^{\pm 1}, \dots, V_{ |\U|}^{\pm 1}  \}$ and, by Theorem 1, $\la  V_1,
\dots, V_{ |\U|+2 }  \ra$ is a $Q$-subgroup of $\G_U( \G_{A,0} ) /  \G_U(
\G_{A,0} )^n$ isomorphic to $B(|\U |+2, n)$. This implies that the subgroup
$\la   V_{ |\U|+1 }, V_{ |\U|+2 } \ra$ of the quotient  $ \G_U( \G_{A,\infty}
) / \G_U( \G_{A,\infty} )^n$ is isomorphic to $B(2,n)$. Observe that any word
in $\U^{\pm 1} \cup t^{\pm 1}$ with 0 sum of exponents on $t$ is equal in
$\HH_t$ to a word of the form $t^{-k} W(\U) t^k$ with $k \ge 0$. Hence, the
kernel $\K_t$ of the epimorphism $\psi_t : \K_t \to \la t \ra_{\infty}$,
given by $\psi_t(t)=t$, $\psi_t(u)=1$, $u \in \U$,  is a group of exponent
$n$ which contains a subgroup isomorphic to $B(2,n)$. Since $V_1, \dots,
V_{|\U|}$ are in the verbal subgroup $\V( \G_{A,0})$  and every element of
$\HH_t$ is conjugate to a word of the form $W(V_1, \dots, V_{ |\U| } )$, it
follows that $\K_t = \V(\K_t)$.

Now we consider the following finite presentation
\begin{equation}
\label{GT2}
\G_{T2} = \la \U , \; t \; \| \; (\ref{R1})-(\ref{R7}), \; t u_j t^{-1} =
V_j, \; j=1, \dots, |\U| \ra
\end{equation}
(recall that  $V_1, V_2, \dots$  are $\A$-words) and observe that, by
definitions and an obvious induction on $i \ge 0$, we have that if $R \in
\R'_{A,i}$ then $R = 1$ in $\G_{T2} $. Hence, $R=1$ in $\G_{T2} $ for every
$R \in \cup_{i=0}^\infty \R'_{A, i}$. By Lemma 1, relations
(\ref{R1})--(\ref{R7}) guarantee that the subgroup $\la V_1, \dots, V_{ |\U|
} \ra$ of $\G_{T2}$ has exponent $n$ and so the subgroup $\la \U \ra$ of
$\G_{T2}$ also has exponent $n$. Now we can see that all defining relations
of the group $\HH_t$  given by (\ref{HHt}) hold in $\G_{T2}$. On the other
hand, it is obvious that all defining relations of $\G_{T2}$ hold in $\HH_t$.
Thus, the groups   $\G_{T2}$ and $\HH_t$ are naturally isomorphic and the
groups  $\G_{T2} = \HH_t$, $\K_{T2} = \K_t$ have all desired properties.
Theorem 2 is proved. \qed
\smallskip

{\em Proof of Theorem 3.} This is similar to the proof of Theorem B
\cite{I02a}. Let $n > 2^{48}$ be odd, $\A(0)$ be an alphabet with $|\A(0)|
>1$ and $G(0) = \la \A(0) \; \| \; R=1, R \in \R(0) \ra$,  where  $\R(0)$
is empty, so that $G(0) = F(\A(0))$ is the free group in $\A(0)$. By
induction on $i \ge 1$, we will construct the alphabet $\A(i)$, the relator
set $\R(i)$,  and the group $G(i) = \la \A(i) \; \| \; R=1, R \in \R(i) \ra$.
To do this, we apply construction of presentation (\ref{GU}) with $\A =
\A(i-1)$, $\R_A = \R(i-1)$ and $\G_A = G(i-1)$.  Then $\A(i) = \U$,
$$
G(i) = \la \A(i) \; \| \; (\ref{R1})-(\ref{R7}), \; R=1, R \in \R_A = \R(i-1)
\ra
$$
and the set $\R(i)$  consists of  all words $R$, where $R=1$ are rewritten in
this form relations of $G(i)$.  By induction and Theorem 1, the groups $G(0)
/ G(0)^n = B(\A(0), n)$, $\dots$, $G(i-1) / G(i-1)^n$ are all naturally embed
in $G(i)$. Therefore, the group
$$
\G_{T3} = G(\infty) =  \la \A(\infty) = \cup_{i=0}^\infty \A(i) \; \| \; R=1,
\; R \in \R(\infty) = \cup_{i=0}^\infty \R(i) \ra
$$
is a group of exponent $n$ such that  all groups
$$
G(0) / G(0)^n = B(\A(0), n),   \dots, \ \  G(i) / G(i)^n, \dots
$$
are naturally embed in $\G_{T3} = G(\infty)$. In particular, $\G_{T3}$
contains a subgroup isomorphic to $B(2,n) \subseteq   B(\A(0), n)$ and if
$W_1, \dots, W_k$ are some $\A(\infty)$-words, say $\A(j)$-words, $j \ge 0$,
then the subgroup $\la   W_1, \dots, W_k \ra $ of $\G_{T3}$ is naturally
isomorphic to the subgroup $\la W_1, \dots, W_k \ra $ of finite
subpresentation
$$
G(j+1) = \la \A(j+1) \; \| \; R=1, R \in \R(j+1) \ra ,
$$  of $\G_{T3}$. Thus, $\G_{T3}$ is a weakly finitely
presented group of exponent $n$ which contains  a subgroup
isomorphic to $B(2,n)$. Theorem 3 is proved. \qed

\end{document}